\documentstyle[preprint,aps]{revtex}
\begin{document}

\draft

\preprint{ SNUTP/96-071}
\title{ On the  $q$-deformed oscillator algebras:
$su_q(1,1)$ and $su_q(2)$}
\author{Phillial Oh\cite{poh}}
\address{Department of Physics\\
Sung Kyun Kwan University\\
Suwon 440-746,  Korea}

\author{Chaiho Rim\cite{rim}}
\address{Department of Physics\\
Chonbuk National University\\
Chonju, 561-756,  Korea}

\date{\today}

\maketitle

\begin{abstract}

We study the relations between $q$-deformations and 
$q$-coherent states of the single oscillator
representations for $su_q(1,1)$ and $su_q(2)$ algebras;
Dyson and  Holstein-Primakoff type in terms of
Biedenharn, Macfarlane and anyonic oscilators.
We also discuss the related Fock-Bargmann
$q$-derivative and integration.
\end{abstract}

\pacs{PACS numbers: 02.90.+p, 03.65.Fd}

\newpage
\def\theequation{\arabic{section}.\arabic{equation}}
\section{Introduction}
\setcounter{equation}{0}

$q$-deformation of oscillator algebra 
was first introduced by Biedenharn (B) \cite{bied} 
and Macfarlane (M) \cite{macf} in the context of
oscillator realization of quantum algebra $su_q(2)$ \cite{skly}.
They obtained their $q$-deformation of $su(2)$ algebra
by using  the double oscillators realization of Jordan-Schwinger 
type. Later, Kulish and Damaskinsky \cite{Kulda} 
used a single oscillator realization of $su(1,1)$ to 
obtain the $q$-deformed algebra
for a special value of Casimir constant equal to $-{3 \over 16}$.
Another oscillator realization of $su(1,1)$ and its 
$q$-deformation was obtained for 
Calogero-type oscillator \cite{cho}.  
This oscillator 
does not satisfy the Heisenberg algebra but a modified one
due to exchange operator, {\it i.e.\/}, Dunkle operator. 
$q$-deformation of this oscillator was also achieved 
and gave a different commutation relation from those of 
B and M.

In this paper, we consider the $q$-deformation of
a single oscillator representation
of $su(1,1)$ and $su(2)$ algebras and their coherent states
in terms of Holstein-Primakoff (HP) \cite{hols}
and Dyson (D) \cite{dyso} realizations which 
appear in many applications from spin density wave in 
condensed matter physics \cite{kitt} to nuclear physics \cite{rowe}.
In spite of the existing literature \cite{bied1}
on the $q$-deformation and $q$-coherent states \cite{gong} of
HP and D realizations of those algebras, 
we find that a comprehensive study on 
the relations among the various realizations is still lacking,
especially the relations among the measures of the $q$-coherent 
state in the resolution of unity. We are going to resolve the issues
by studying other types of $q$-deformations suited to HP and D
realizations, starting with the "symmetric" $q$-deformation of the
Lie algebras developed by Curtright and Zachos \cite{curt}.
We find that this procedure is instructive since we can obtain the 
B, M, $q$-deformed anyonic 
oscillators, and their $q$-coherent states
naturally. Also, $q$-deformation of Fock-Bargmann (FB) type 
can be realized in $q$-derivatives.

This paper is organized as follows. 
In section II, we obtain $su_q(1,1)$ in HP and D 
realizations in terms of the B, M,
and anyonic type oscillators.
We construct the $q$-coherent states 
and compare the measure in each case.
We also present the FB representation in $q$-derivative. 
In Section III, the same analysis is performed for $su_q(2)$ case.
Section IV contains the conclusion and discussion.

Before going into details, we briefly explain our notation.
In the oscillator representation of the Heisenberg-Weyl algebra,
\begin{equation}
[a_-, a_+] =1\,,
\label{heisenberg}
\end{equation}
we introduce a number eigenstate $|n>$ 
of number operator $N=a_+ a_-$.
We require 
$|0>$ be  annihilated by $a_-$,
$a_- |0> = 0$.
Explicitly, the creation and annihilation operators 
act on the ket, 
\begin{equation}
a_-|n> = n |n-1>\,
\quad
a_+ |n> = |n+1>\,.
\label{anrelation}
\end{equation}
It should be noted that the normalization of the ket 
is not fixed yet, which  will result
 from Hermitian property of the generators of 
 $su_q(1,1)$ and $su_q(2)$. We use 
FB holomorphic representation of the 
oscillator algebra, in which
$ < \xi | n> = \xi^n$ and 
\begin{equation}
{d \over d \xi} < \xi| n> = < \xi | a_- | n> \,
\quad 
\xi < \xi| n> = < \xi | a_+ | n>  \, .
\end{equation}

\def\theequation{\arabic{section}.\arabic{equation}}
\section{$su_q(1,1)$ and coherent state}
\setcounter{equation}{0}

$su(1,1)$ satisfies the algebra,
\begin{equation}
[ K_0,  K_{\pm}]=\pm K_{\pm}, \quad 
[ K_+,  K_-]=-2K_0,
\label{su11algebr}
\end{equation}
and Casimir invariant is expressed as
$ C=K_0(K_0-1)-K_+K_-\,$.
To have a connection with the  oscillator algebra,
we require the number eigenstate $|n>$ be
an eigenstate of $K_0$,
\begin{equation}
K_0 |n> = (k_0 + n) |n>\,.
\end{equation}
Here, $k_0$ is assumed to be
a positive integer or half odd integer
and $|0>$  is also annihilated by $K_-$.
This representation gives the Casimir number
$k_0(k_0 -1)$.
For definiteness, we will assume that 
the $su(1,1)$ algebra act on the ket $|n>$ as 
\begin{equation}
K_- |n> = n | n-1>\,\quad
K_+ |n> = (n + 2 k_0) |n + 1>\,.
\label{k-ladder}
\end{equation}
This convention is consistent with 
the holomorphic first-order differential operator 
representation for $su(1,1)$ given in \cite{pere} 
\begin{equation}
\hat K_+ (\xi) =  \xi^2 {d \over d \xi} + 2 k_0  \xi\,,\quad
\hat K_-(\xi) = {d \over d \xi} \,,\quad
\hat K_0 (\xi) = \xi {d \over d \xi} + k_0 \,.
\label{su11holo}
\end{equation}
Since $
\hat K_i(\xi) < \xi| n> = < \xi | K_i  | n>
$, we can check the relation in Eq.~(\ref{k-ladder}) holds. 

$q$-deformed algebra $su_q(1,1)$ is given as \cite{skly}
\begin{equation}
 [Q_0, Q_{\pm}]=\pm  Q_{\pm}\,,\quad
[ Q_+,  Q_-]=-[2Q_0]_{q}\,,
\label{su11algebra}
\end{equation}
where the $q$-deformation is defined as
\begin{equation}
[x]_q\equiv \frac{q^x-q^{-x}}{q-q^{-1}}\,.
\end{equation}
$q$-deformed Casimir invariant is given by 
$C_q = [Q_0]_q[Q_0 - 1]_q - Q_+ Q_-$.
One can obtain the explicit form of the $q$-deformed 
generators following \cite{curt},
\begin{equation}
Q_0=K_0\,,\quad 
Q_-=K_-f(K_0)\,,\quad
Q_+=f(K_0)K_+\,.
\label{su11qrealization}
\end{equation}
Noting a useful identity
$ -[2K_0]_{q}=g(K_0)-g(K_0+1)$
where 
$g(K_0)=[K_0-k_0]_{q}[K_0+k_0-1]_{q}\,$
and
$g(K_0=k_0)=0\,$,
we  may identify $g(K_0)$ as $f(K_0)^2 (K_0(K_0-1)-C)$,
\begin{equation}
f(K_0) = \sqrt{[K-k_0]_q [K_0 + k_0 -1]_q
	\over (K-k_0) (K + k_0 -1)}
	= \sqrt{[N]_q [N + 2k_0 -1]_q
	\over N (N + 2k_0 -1)}\,.
\label{fk0}
\end{equation}

This realization is not the unique choice.
In the following, we give some other examples
which will result in the  $q$-deformed oscillator algebra
or the $q$-deformed FB representation.
In addition, independent of the realization
we are requiring the conjugate
relations,
\begin{equation}
Q_-^{\dagger} = Q_+\,,\quad
Q_0^{\dagger} = Q_0\,.
\label{conjugate}
\end{equation}
This will determine the norm of each state
for the given realization and the resolution of 
unity for coherent state.

\subsection{B oscillator realization.}

(1) D type.

Let us consider the realization,
\begin{equation}
Q_0=K_0 = N+ k_0\,,\quad 
Q_-=K_- \sqrt{[N]_q \over N} \,,\quad
Q_+=\sqrt{[N]_q \over N} 
{[N + 2k_0 -1]_q \over (N+ 2k_0 -1)} K_+\,.
\label{su11-ID}
\end{equation}
This is obtained if we re-scale 
$Q_+$ and $Q_-$ in Eq.~(\ref{su11qrealization}).
This ladder  operators act on the ket as 
\begin{equation}
Q_- |n>_D = \sqrt{n [n]_q} |n-1>_D \,,\qquad
Q_+ |n>_D = \sqrt{ [n+1]_q \over n+1} \, [n+ 2k_0]_q\, |n+1>_D\,.
\end{equation}
Eq.~(\ref{su11-ID}) becomes 
an oscillator realization if we interpret it as
\begin{equation}
Q_0 = N+ k_0\,,\quad
Q_- = (a_q)_- \,,\quad
Q_+ = [N + 2k_0-1 ]_q (a_q)_+ \,,
\end{equation}
and identify $(a_q)_\pm$ as 
\begin{equation}
(a_q)_- = a_-  \sqrt{[N]_q \over N}\,,\quad
(a_q)_+ =  \sqrt{[N]_q \over N}\,\, a_+ \,.
\label{a-Bied}
\end{equation}
$(a_q)_\pm$ satisfy  the $q$-deformed 
oscillator algebra of B type,
\begin{equation}
(a_q)_- (a_q)_+ - q (a_q)_+ (a_q)_- = q^{-N}\,.
\label{Biedenharn}
\end{equation}

$q$-deformed  coherent state is defined  by 
( {\it \`a la} Perelomov' coherent 
state \cite{pere})
\begin{equation}
|z>_D = e_q^{\bar z Q_+} |0>_D
=\sum_{n=0}^\infty \bar z^n 
\sqrt{1 \over [n]_q! n!} 
{[n+2k_0-1]_q! \over [2k_0 -1]_q! } \,\, |n>_D\,,
\label{su11qcoherent}
\end{equation}
where we use the $q$-deformed exponential function.
The subscript D stands for D type.
(Note that we do not add a  normalization constant 
in this definition since this will introduce $z$ in addition to 
$\bar z$). 
Conjugate relation, Eq.~(\ref{conjugate}) gives the normalization 
of the number eigenstate,
\begin{equation}
{}_D\!\!<n|n>_D = {n! [2k_0 -1]_q! \over [n+2k_0 -1]_q!}\,.
\label{normalization}
\end{equation}
This normalization provides us the resolution of unity as
\begin{equation}
I = \sum_{n=0}^{\infty}  
{[n + 2k_0 -1]_q! \over  n! [ 2k_0 -1]_q!} 
|n>_D\,{}_D\!\!<n|
= \int d^2_q z \, G(z)\, |z>_D\,{}_D\!\!<z|\,,
\end{equation}
and the measure  $G(z)$ is given as 
\begin{equation}
G(z)=\left\{ \begin{array}{ll}
{[2k_0 -1]_q \over \pi} (1 - |z|^2)_q^{2k_0 -2}&
 \text{for } 2k_0= \text{integer} >1 \\
 {|z|^2 \over \pi} & \text{for } k_0= 1
 \end{array}\right.
\label{su11qmeasure}
\end{equation}
where the $q$-deformed function is defined as 
$(1-x)^n _q = \sum_{m=0}^n {[n]_q! \over [m]_q! [n-m]_q!} (-x)^m$.
(For $k_0= {1 \over 2}$, see below Eq. (\ref{aq-measure})).
Here,  the two dimensional integration is defined as 
\begin{equation}
d^2_q z \equiv {1 \over 2} d\theta\,\, d_q |z|^2 \,.
\label{integ}
\end{equation}
The angular integration is an 
ordinary integration, $0 \le \theta \le 2\pi$.
The radial part is a  $q$-integration, which is the inverse 
operation of $q$-derivative defined as
\begin{equation}
{d \over d_q z} f(z) = {f(qz) - f(q^{-1}z) \over z(q- q^{-1})}\,.
\label{qderivative}
\end{equation}
One can check that $I$ commutes with the $su_q(1,1)$ generators.

We note that in this Hilbert space, 
$(a_q)_+$ is not an adjoint of $(a_q)_-$.
To have the conjugation property
between $(a_q)_-$ and $(a_q)_+$ 
as well as between $Q_-$ and $Q_+$,
we may resort to HP realization.
Therefore, we need to compare the quantities in different 
realization. 
It turns out to be useful 
to express the quantities in 
unit normalized eigenstate  basis $|n)$ instead of 
$|n>$. 
For later comparison, we give the explicit expression;
\begin{equation}
Q_- |n)_D = \sqrt{[n]_q [n+2 k_0 -1]_q} |n-1)_D\,,\quad
Q_+ |n)_D = \sqrt{[n+1]_q [n+2 k_0 ]_q} |n+1)_D\,.
\label{Qrelation}
\end{equation}
And the coherent state in Eq.~(\ref{su11qcoherent}) becomes
\begin{equation}
|z>_D = e_q^{\bar z Q_+} |0>_D
=\sum_{n=0}^\infty \bar z^n 
\sqrt{ [n+2k_0-1]_q! \over [n]_q! [2k_0 -1]_q! } \,\, |n)_D\,.
\label{su11qcoherent-0}
\end{equation}

(2) HP type.

Let us consider the realization in terms of $q$-deformed oscillator 
in Eq.~(\ref{a-Bied}) to  have the HP realization,
\begin{eqnarray}
&&Q_0= N+ k_0\,,
\nonumber \\
&& Q_-=K_- \sqrt{[N]_q [N+2k_0-1]_q\over N }
   = (a_q)_- \sqrt{[N +2k_0 -1]_q} \,,
\nonumber \\
&&Q_+=\sqrt{[N]_q [N + 2k_0 -1]_q \over N} 
       {1 \over (N+ 2k_0 -1)} K_+
   = \sqrt{[N + 2k_0 -1]_q} (a_q)_+\,.
\label{su11-IHP}
\end{eqnarray}
The ladder operators act on the ket as
\begin{eqnarray}
Q_- |n>_H &=& 
\sqrt{ n [n]_q [n+2k_0 -1]_q } |n-1>_H \nonumber\\
Q_+ |n>_H &=& \sqrt{ [n+1]_q [n+2k_0]_q \over n+1} 
	  \, |n+1>_H.
\label{QHP-relation}
\end{eqnarray}
The subscript $H$ stands for HP.

Conjugate relation between $Q_\pm$'s  requires 
the normalization in this 
Hilbert space,\\
${}_H\!\!<n|n>_H = n!$,
which is different from the previous one,
Eq.~(\ref{normalization}).
This  is not surprising since 
two Hilbert spaces are different.
If we use the unit normalized ket $|n)_H$,  
then the relation given in Eq.~(\ref{QHP-relation}) becomes 
exactly the same form given in Eq.~(\ref{Qrelation})
with subscript D replaced by subscript H, and the $q$-deformed 
coherent state corresponding  
to Eq.~(\ref{su11qcoherent}) is given as  
\begin{equation}
|z>_H = e_q^{\bar z Q_+} |0>_H
=\sum_{n=0}^\infty \bar z^n 
\sqrt{[n + 2k_0 -1]_q! \over [n]_q! [2k_0 -1]_q!} 
\,\, |n)_H
\end{equation}
in terms of the normalized ket $|n)_H$. 
This is again exactly the same form given 
in Eq.~(\ref{su11qcoherent-0}).
Therefore, the resolution of unity for the coherent state 
is expressed  in terms of the same measure $G(z)$ 
in  Eq.~(\ref{su11qmeasure}),
even though two realizations look so different at first sight.

Note that in this HP realization, 
the conjugate relation 
between $a_+^\dagger = a_-$ is satisfied automatically,
since 
\begin{equation}
(a_q)_- |n)_H = \sqrt{[n]_q}\,|n-1)_H\,\quad
(a_q)_+ |n)_H = \sqrt{[n+1]_q} |n+1)_H \,.
\label{aq-relation}
\end{equation}
Therefore, one can equally use a new coherent state,
$q$-deformed version of Glauber coherent state \cite{klau},
\begin{equation}
|z]_H = e^{\bar z a_+}_q |0>_H\,.
\label{qGlaubercoherent}
\end{equation}
Explicit form of this coherent state is given as
$ |z]_H = \sum_{n=0}^{\infty} {\bar z^n \over [n]_q!} |n)_H\,$.
The resolution of unity is given as
\begin{equation}
I = \sum_{n=0} ^\infty  |n)_H\,{}_H\!(n|
= \int d^2_q z\,  g(z) | z]_H\, {}_H\![ z|\, 
\end{equation}
where $g(z)$ is given by 
\begin{equation}
g(z) = {1 \over \pi} e_q^{- |z|^2}\,,
\label{aq-measure}
\end{equation}
and the domain is over an infinite plane. 
Unlike the D case, 
this holds for any value of $k_0$.
Therefore, 
for $k_0={1 \over 2}$, 
one can use the $q$-analogue of Glauber 
coherent state. 
For other value of $k_0$, 
one can use the Bargmann measure defined in 
Eq.~(\ref{aq-measure}) or 
Liouville measure in 
Eq.~(\ref{su11qmeasure})  
depending on the definition of coherent state.

\subsection{M oscillator realization.}

(1) D type.

We  consider  a different realization from the B
type oscillator realization,
\begin{eqnarray}
&&Q_0=K_0= N + k_0 \,,
\nonumber\\
&&Q_-=K_- \sqrt{[N]_q \over N}q^{N-2 \over 2}= (b_q)_-\,,
\nonumber\\
&&Q_+=q^{-{N -1 \over 2}} \sqrt{[N]_q \over N} 
{[N + 2k_0 -1]_q \over (N+ 2k_0 -1)} K_+
=q^{-(N-1)} [N + 2k_0 -1]_q (b_q)_+ \,,
\label{su11-IID}
\end{eqnarray}
with new oscillator given as 
\begin{equation}
(b_q)_- = (a_q)_- q^{N-1 \over 2}
	= a_- \sqrt{\{N\}_q \over N}\,,\quad
(b_q)_+ = q^{N-1 \over 2}(a_q)_+
	= \sqrt{\{N\}_q \over N} a_+ \,,
\label{aq-Mac}
\end{equation}
where we introduce a new definition of $q$-number,
\begin{equation}
\{x\}_q = {q^{2x} - 1 \over q^2 -1} = [x]_q \,\, q^{x-1}\,.
\end{equation}
This oscillator realization gives the $q$-deformed 
oscillator algebra of M type,
\begin{equation}
(b_q)_- (b_q)_+ - q^2 (b_q)_+ (b_q)_- = 1\,.
\label{Macfarlane}
\end{equation}

In terms of this realization,  the ladder operators act
on the ket as 
\begin{equation}
Q_- |n> =  \sqrt{n \{n\}_q} |n-1> \,,\qquad
Q_+ |n> = q^{-n }\sqrt{ \{n+1\}_q \over n+1} \, 
[n+ 2k_0]_q\, |n+1>\,.
\end{equation}
In the following, for notational simplicity, 
we will  delete the 
subscript on the ket which distinguishes the Hilbert space,
since there is no possibility of confusion.
The conjugate relation between $Q_\pm$'s 
gives the normalization, 
\begin{equation}
<n|n> = { n! [ 2k_0 -1]_q! \over [n + 2k_0 -1]_q!}
q^{n(n-1) \over 2} \,.
\end{equation}
In terms of unit normalized  ket $|n)$,
we have the canonical operator relations
for $Q_\pm$  as in Eq.~(\ref{Qrelation}).
In addition, $q$-deformed  coherent state 
for $su_q(1,1)$ has the same form as in Eq.~(\ref{su11qcoherent-0})
and therefore, the measure $G(z)$ 
in Eq.(\ref{su11qmeasure})
is used for the resolution of unity 
for the coherent state.

(2) HP type.

Another HP type realization is given as
\begin{eqnarray}
&&Q_0= N+ k_0\,,
\nonumber \\
&&Q_- = (b_q)_- \sqrt{q^{-(N-1)} [N +2k_0 -1]_q} \,,
\nonumber \\
&&Q_+=\sqrt{q^{-(N-1)} [N + 2k_0 -1]_q} (b_q)_+\,,
\label{su11-IIHP}
\end{eqnarray}
with $b_q$'s defined in Eq.~(\ref{aq-Mac}).
However, these generators coincide with the one given in 
HP type of the B oscillator realization,
Eq.~(\ref{su11-IHP}).
That is, the Hilbert space is exactly same for both cases
as far as $su_q(1,1)$ is concerned.
Therefore, the $su_q(1,1)$ $q$-deformed  coherent state 
in terms of the normalized ket $|n)$
is exactly the same form given in Eq.~(\ref{su11qcoherent-0})
and the same measure $G(z)$ 
in Eq.~(\ref{su11qmeasure})
is used for the resolution of unity.

On the other hand, from the oscillator point of view,
one can define a new coherent state since
the conjugate relation 
between $(b_q)_+^\dagger = (b_q)_-$ is satisfied automatically;
\begin{equation}
(b_q)_- |n) = \sqrt{\{n\}_q}\,|n-1)\,\quad
(b_q)_+ |n) = \sqrt{\{n+1\}_q} |n+1) \,.
\label{dagger}
\end{equation}
Let us  define  
another version of $q$-deformed Glauber coherent state as
\begin{equation}
|z\} = E_q^{\bar z (b_q)_+}|0> 
     = \sum_{n=0}^{\infty} {\bar z^n \over \{n\}_q!} |n)\,,
\label{newqcoherent}
\end{equation}
where $q$-deformed exponential function $E_q^x$ differs 
from $e_q^x$ in that $q$-number $[n]_q$ is replaced by
$\{n\}_q$,
\begin{equation}
E_q^x = \sum_{n=0}^\infty {x^n \over \{n\}_q!}
\end{equation} 
The resolution of unity is given as
\begin{equation}
I 
= \sum_{n=0} ^\infty  |n)\,(n|
= \int d^2_q z\,  h(z) | z\}\, \{ z|\,. 
\end{equation}
where $h(z)$ is given by 
\begin{equation}
h(z) = {1 \over \pi} E_q^{- |z|^2}\,,
\label{bq-measure}
\end{equation}
and the domain is over  an infinite plane.

\subsection{$q$-anyonic oscillator realization}

Let us consider the HP type realization again.
As seen in the previous section, one can have the 
B oscillator Eq.~(\ref{su11-IHP}) 
or M oscillator Eq.~(\ref{su11-IIHP})
from the same $q$-deformed form of the $Q_i$'s.
We give another useful form of oscillator realization, 
$q$-anyonic oscillator.
Since D and HP type realizations are 
now trivially connected, we present only HP
realization which maintains the conjugate condition for 
oscillator algebra also.
Let us put $Q_i$'s as 
\begin{eqnarray}
&&Q_0= N+ k_0\,,
\nonumber \\
&& Q_-=K_- \sqrt{[N]_q [N+2k_0-1]_q\over N }
   = (A_q)_- \sqrt{ (A_q)_+ (A_q)_+ + 2[k_0 - {1 \over 2}]_q} \,,
\nonumber \\
&&Q_+=\sqrt{[N]_q [N + 2k_0 -1]_q \over N} 
       {1 \over (N+ 2k_0 -1)} K_+
   = \sqrt{ (A_q)_+ (A_q)_+ + 2[k_0 - {1 \over 2}]_q}\,\, (A_q)_+ \,,
\end{eqnarray}
Then we have the $q$-deformed oscillator as 
\begin{equation}
(A_q)_-= a_- \,\, \sqrt{[N+ k_0 - {1\over 2} ]_q - 
			[k_0 - {1 \over 2}]_q 
		  \over N}\,,
\quad
(A_q)_+= \sqrt{[N+ k_0 - {1\over 2} ]_q - [k_0 - {1 \over 2}]_q 
		  \over N}\,\, a_+\,.
\end{equation}

Its commutation relation looks complicated,
\begin{equation}
((A_q)_- (A_q)_+ + [k_0 - {1 \over 2}]_q)
-q ((A_q)_+ (A_q)_- + [k_0 - {1 \over 2}]_q)
= q^{-(N + k_0 + { 1\over 2})}\,.
\end{equation}
However, the meaning of this commutation relation becomes clear if 
we rewrite the relation in M's form,
\begin{equation}
(B_q)_- (B_q)_+ - q^2 (B_q)_+ (B_q)_- =1 \,,
\label{paraMac}
\end{equation}
by identifying 
\begin{eqnarray}
(B_q)_- (B_q)_+ = q^{N + k_0 - {1 \over 2}}[N +k_0 + {1 \over 2}]_q
	        = q^{N + k_0 - {1 \over 2}}
		\,\,((A_q)_- (A_q)_+ + [k_0 - {1 \over 2}]_q)\,,
\nonumber\\
(B_q)_+ (B_q)_- = q^{N + k_0 - {3 \over 2}}[N+ k_0 - {1 \over 2}]_q
	        = q^{N + k_0 - {3 \over 2}}
		\,\,((A_q)_+ (A_q)_- + [k_0 - {1 \over 2}]_q)\,.
\label{Bquad}
\end{eqnarray}
In this realization, 
the vacuum $|0>$ is not annihilated by $(B_q)_-$ 
unless $k_0 = {1 \over2}$, since 
\begin{equation}
(B_q)_+ (B_q)_- |0> = \{k_0 - {1 \over 2}\}_q |0>\,.
\end{equation}
This feature reflects the fact that this realization corresponds 
to the $q$-deformed non-trivial one dimensional analogue of 
anyon which appears in two dimensional oscillator representation
with $k_0$ being related with statistical
parameter in anyon physics \cite{chorim}. 

The conjugate relation between $(B_q)_-$ and $(B_q)_+$
can be seen  formally at the operator level in Eq.~(\ref{Bquad})
since the conjugate relation between $(A_q)_-$ and $(A_q)_+$ does
hold. However, the fact that $(B_q)_-$ does not annihilate the 
vacuum $|0>$ implies that one cannot define a proper Hilbert space.
Therefore, the measure of the coherent state 
of the Glauber type for the $B_q$ oscillator system cannot be defined.
On the other hand, the measure of the $q$-deformed 
coherent state of Perelomov type is given in Eq.~(\ref{su11qmeasure}).
One may also define the coherent state of the 
Glauber type in terms of 
the $A_q$ oscillator, whose explicit form of the measure 
turns out to be very complicated and will not be reproduced here.

As far as $B_q$ oscillator is concerned, we may 
construct a Hilbert space from a new vacuum which is 
annihilated by $(B_q)_-$. Then, since the commutation 
relation Eq.~(\ref{paraMac}) is the same form as  
in Eq.~(\ref{Macfarlane}),
the generators act on the 
new Hilbert space as in Eq.~(\ref{dagger}).
In this case, one can contruct the $q$-deformed Glauber type 
coherent state and the measure is given in Eq.~(\ref{bq-measure}).
However, this representation has  nothing to do with the 
anyonic representation mentioned above.

We comment in passing that there is another well-known 
one dimensional oscillator representation for anyon type;
Calogero oscillator system, which turns out to be 
the realization \cite{cho,mac} of parabose system \cite{gree}. 
Its $q$-deformed realization does not 
satisfy the commutation relation of
M type Eq.~(\ref{paraMac}). 
The explicit measure for the $q$-deformed 
coherent state of the Glauber type in this case is
already known \cite{cho}.

\subsection{FB realization with symmetric $q$-derivative}

Let us consider a realization,
\begin{equation}
Q_0=K_0 = N+ k_0\,,\quad 
Q_-=K_- {[N]_q \over N} \,,\quad
Q_+= {[N + 2k_0 -1]_q \over (N+ 2k_0 -1)} K_+\,.
\label{su11-IIIB}
\end{equation}
These generators act on the ket as 
\begin{equation}
Q_- |n> =  [n]_q |n-1> \,,\qquad
Q_+ |n> = [n+ 2k_0]_q\, |n+1>\,.
\end{equation}
Conjugate relation between $Q_\pm$'s
requires the normalization 
of the number eigenstate,
\begin{equation}
<n|n> = {[n]_q! [2k_0 -1]_q! \over [n+2k_0 -1]_q!}.
\end{equation}
In terms of the unit normalized ket $|n)$,
we reproduce the same form of 
$su_q(1,1)$ coherent state and 
resolution of unity as seen in the previous subsections, A and B.

What makes this  realization  different from the previous ones 
is that it gives a  natural $q$-deformation of the FB representation
of $su(1,1)$. By using  $<\xi|n>= \xi^n$, we have
\begin{equation}
\hat Q_+ (\xi) =  \xi [\xi {d \over d \xi} + 2 k_0 ]_q \,,\quad
\hat Q_-(\xi) = {d \over d_q \xi} \,,\quad
\hat Q_0 (\xi) = \xi {d \over d \xi} + k_0 \,.
\end{equation}
The $q$-derivative in $\hat Q_-$ is 
defined in Eq. (\ref{qderivative}).
This implies that the oscillator realization is given as
\begin{equation}
(a_q)_-(\xi) = {d \over d_q \xi} \,,\quad
(a_q)_+(\xi) = \xi \,,
\label{axi}
\end{equation}
which satisfies the $q$-deformed oscillator algebra 
of B type, Eq.~(\ref{Biedenharn}).

\subsection{FB realization with a-symmetric $q$-derivative}

We may consider  a little  modified version of Eq.~(\ref{su11-IIIB}),
\begin{equation}
Q_0=K_0 = N+ k_0\,,\quad 
Q_-=K_- {[N]_q \over N}q^{N-1} \,,\quad
Q_+= q^{-(N-1)}{[N + 2k_0 -1]_q \over (N+ 2k_0 -1)} K_+\,.
\label{su11-IIIFB}
\end{equation}
Then the generators act on the ket as 
\begin{equation}
Q_- |n> =  \{n\}_q |n-1> \,,\qquad
Q_+ |n> = q^{-n} [n+ 2k_0]_q\, |n+1>\,.
\end{equation}
Conjugate relation between $Q_\pm$'s 
requires the normalization 
of the ket as,
\begin{equation}
<n|n> = q^{n(n+2k_0-2)}
	 {\{n\}_q! \{2k_0 -1\}_q! \over \{n+2k_0 -1\}_q!}
\end{equation}
One can easily check that in this Hilbert space,
the same form of $su_q(1,1)$  coherent state and 
resolution of unity are reproduced 
as in the previous sections
if we use the unit normalized ket $|n)$.

We have a similar  
$q$-deformation of the FB representation
of $su(1,1)$ as in the previous section,
\begin{equation}
\hat Q_0 (\xi) = \xi {d \over d \xi} + k_0 \,,\quad
\hat Q_-(\xi) = {D \over D_q \xi} \,,\quad
\hat Q_+ (\xi) =  \xi q^{-(2 \xi {d \over d \xi} + 2k_0 -1)}
	\{\xi {d \over d \xi} + 2 k_0 \}_q \,.
\end{equation}
The derivative in $\hat Q_-$ is replaced by a new  $q$-derivative
which is given by 
\begin{equation}
{D \over D_q z} f(z) = {f(q^2z) - f(z) \over z(q^2- 1)}\,,
\end{equation}
This implies that the oscillator realization is given by
\begin{equation}
(b_q)_-(\xi) = {D \over D_q \xi} \,,\quad
(b_q)_+(\xi) = \xi \,,
\end{equation}
which satisfies the $q$-deformed oscillator algebra 
of M type, Eq.~(\ref{Macfarlane}).

\def\theequation{\arabic{section}.\arabic{equation}}
\section{$su_q(2)$ and coherent state}
\setcounter{equation}{0}

$su_q(2)$ and its coherent state can be studied 
in close analogy with 
the previous section and therefore,
we will describe briefly about B oscillator 
realization only.
$su(2)$ satisfies the algebra,
\begin{equation}
[ K_3,  K_{\pm}]=\pm K_{\pm}, \quad 
[ K_+,  K_-]=2K_3,
\label{su2algebra}
\end{equation}
and Casimir operator is given  as
$C=K_3(K_3+1)+K_-K_+.$ In addition, 
\begin{equation}
K_- |n> = n |n-1>\,, \quad
K_+ |n> = (J-n) |n+1>\,.
\end{equation}
The Hilbert space is finite dimensional with dimension
$J +1 $ where 
\begin{equation}
K_-\vert n=0>=0\,,\quad
K_+ \vert n=J> = 0\,,
\end{equation}
where $J$ is an integer.
Since $|n>$ is an eigenstate of $K_3$,
\begin{equation}
K_3 |n>= (n - {J \over 2})|n>\,,
\end{equation}
we have the Casimir constant,
$C= {J \over 2}({J \over 2} + 1)$.

$q$-deformed $su(2)$ algebra is given as \cite{skly}
\begin{equation}
[ Q_3,  Q_{\pm}]=\pm  Q_{\pm}\,,\quad 
[Q_+, Q_-]=[2 Q_3]_{q^2}\,.
\end{equation}
We require  conjugate relation $Q_-^{\dagger} = Q_+$,
$Q_3^{\dagger} = Q_3$  
independent of the realization.
Repeating the same procedure as $su(1,1)$ case, we find
\begin{equation}
Q_3= K_3\,,\quad
Q_-= K_- F(K_3)\,,\quad
Q_+ = F(K_3) K_+\,.
\label{su2qcom}
\end{equation}
where
\begin{equation}
F(K_3)=
\sqrt{[{J \over 2} + K_3]_q [{J \over 2} +1 - K_3]_q
	\over ({J \over 2} + K_3)({J \over 2} +1 - K_3)}
= \sqrt{[N]_q [J +1 - N]_q \over N(J +1 - N)}\,.
\end{equation}

\subsection{D type of B oscillator representation.}

\begin{equation}
Q_3= K_3 = {J \over 2} - N\,,\quad
Q_-= K_- \sqrt{[N]_q \over N}\,,\quad
Q_+ = \sqrt{[N]_q \over N} {[J+1-N]_q \over (J+1-N)} K_+\,.
\label{su2-D}
\end{equation}
These generators act on the ket as 
\begin{equation}
Q_- |n> = \sqrt{n [n]_q} |n-1> \,,\qquad
Q_+ |n> = \sqrt{ [n+1]_q \over n+1} \, [J-n]_q\, |n+1>\,.
\end{equation}
We get an oscillator realization if 
\begin{equation}
Q_0 = {J \over 2} - N\,,\quad
Q_- = (a_q)_- \,,\quad
Q_+ = [J - N+1   ]_q (a_q)_+ \,.
\end{equation}
$(a_q)_\pm$ is as defined in Eq.~(\ref{a-Bied}).

Introducing  unit normalized ket,
$ |n) = \sqrt{J! \over n1 (J-n)!}\,|n>\,$,
we have the canonical operator relations of $su_q(2)$.
\begin{equation}
Q_- |n) = \sqrt{[n]_q [J+1-n ]_q} |n-1)\,,\quad
Q_+ |n) = \sqrt{[n+1]_q [J-n ]_q} |n+1)\,.
\label{su2Qrelation}
\end{equation}
$q$-deformed  coherent state is given as 
\begin{equation}
|z> = e_q^{\bar z Q_-} |J>
=\sum_{n=0}^\infty \bar z^n 
\sqrt{[J]_q! \over [n]_q! [J-n]_q!}\,|n)\,. 
\label{su2qcoherent}
\end{equation}
Resolution of unity is expressed as
\begin{equation}
I = \sum_{n=0}^{J}  
|n)\,(n|
= \int d^2_q z \, H(z)\, |z>\,<z|\,,
\end{equation}
and the measure is given as 
\begin{equation}
H(z)= {[J+1]_q \over \pi} {1 \over (1 + |z|^2)_q^{2+J}}\,.
\label{su2qmeasure}
\end{equation}
One can check that $I$ commutes with the $su_q(2)$ generators.

\subsection{HP 
type of B oscillator representation.}

\begin{eqnarray}
&&Q_3= K_3 = {J \over 2} - N\,,
\nonumber \\
&&Q_-= K_- \sqrt{[N]_q [J+1-N]_q \over N}
   = (a_q)_-\sqrt{[J+1-N]_q}\,,
\nonumber \\
&&Q_+ = \sqrt{[N]_q [J+1-N]_q \over N} 
    {1 \over (J+1-N)}  K_+
    = \sqrt{[J+1-N]_q} (a_q)_+ \,.
\label{su2-HP}
\end{eqnarray}
$(a_q)_\pm$ is defined in Eq.~(\ref{a-Bied}).
The ladder operators act on the ket as 
\begin{equation}
Q_- |n> = \sqrt{n [n]_q [J +1 -n]_q} |n-1> \,,\qquad
Q_+ |n> = \sqrt{ [n+1]_q [J-n]_q \over n+1} \, |n+1>\,.
\end{equation}
Using the normalized ket, 
$|n) = \sqrt{1\over n! }\,|n>\,$,
we have the canonical 
ladder operator realization as in Eq.~(\ref{su2Qrelation})
and  $q$-deformed coherent state is given as 
\begin{equation}
|z> = e_q^{\bar z Q_-} |J>
=\sum_{n=0}^\infty \bar z^n 
\sqrt{[J]_q! \over [n]_q! [J-n]_q!}\,|n)\,. 
\end{equation}
Therefore, the measure $H(z)$ 
given in Eq.(\ref{su2qmeasure})
is used for the resolution of unity.

Because of the conjugate relation 
between $(a_q)_+$ and $(a_q)_-$,
we can equally consider the $q$-coherent state 
of finite Glauber coherent state.
However, the Hilbert space is finite dimensional, so
one has to modify the definition of the coherent
state from the $su(1,1)$ case, 
Eq.~(\ref{newqcoherent});
\begin{equation}
|z>> = e^{\bar z (a_q)_-}_q |J)
= \sum_{n=0}^{\infty} \bar z^n 
 \sqrt{[J]_q! \over [n]_q! [J-n]_q!}\,\, |n)\,.
\end{equation} 
It is interesting to note that the 
oscillator coherent state 
reproduces the same form of $su_q(2)$ coherent state 
given in Eq.~(\ref{su2qcoherent}).
This is because the Hilbert space is finite dimensional 
 in contrast with $su_q(1,1)$ case.
\def\theequation{\arabic{section}.\arabic{equation}}
\section{conclusion}
\setcounter{equation}{0}

We have presented and compared various 
type of oscillator algebra  realizations of 
$su_q(1,1)$  and $su_q(2)$ algebras, and their coherent states.
For $su_q(1,1)$, if  we impose the conjugate condition 
for the generators, the Perelomov $q$-coherent states has 
a common measure in the resolution of unity independently
of the explicit forms of realization. Another type of
$q$-coherent state, the Glauber type is considered in the HP
realization since $a_-$ and $a_+$ are automatically conjugate
to each other. The explicit measure for this type of $q$-coherent
state depends on the oscillator realization such as
B or M type; the Liouville type measure
 defined in  Eq.~(\ref{su11qmeasure}) or  
the Bargmann type in Eq.~(\ref{aq-measure}),
or the other Bargmann type in Eq.~(\ref{bq-measure}).
 In addition, it is shown that
the explicit forms of the generators of $su_q(1,1)$ can be
modified such that $q$-anyonic oscillator and 
various definition of $q$-derivatives can be 
accommodated in the realizations. 

$su_q(2)$ shares much of the same results with $su_q(1,1)$.
However,  in HP realization, the finite Glauber $q$-coherent 
state does not have the Bargmann measure, but has the 
Liouville measure. The difference comes from the finiteness
of the dimension of the Hilbert space. Therefore, the measure
of coherent state in $su_q(2)$ is distinguishable from that
of the oscillator coherent state on a plane.

We also presented  two different types of FB realization which
provide two different definitions of $q$-derivative and 
$q$-integration such that we can describe their $q$-deformed 
oscillators algebra in a natural and simple fashion.

We conclude with a couple of  remarks. 
The D representations  can be extended 
to the $SU(N)$ case \cite{oh955}. 
The HP version in the
$SU(N)$ case can also be constructed \cite{rand}. 
In addition, its $q$-deformation was considered in  \cite{sunfu}.
It would be interesting to go through the same analysis in this
higher case, especially in connection with FB realization.
 
$q$-deformed FB representation will be useful for evaluating
the $q$-deformed version of the path integral \cite{baul}. 
In our approach, $q$-deformation of FB representation
is understood in terms of the oscillator representation
and the role of the $q$-derivatives are illustrated.
However, $q$-integration is performed essentially 
for one dimensional direction, radial part.
Angular part is treated as an ordinary integration 
Eq. (\ref{integ}). So our
resolution of unity cannot be used directly in evaluating
the $q$-deformed version of path integral at this stage.
To overcome the shortcomings, one has to fully develop $q$-deformed 
higher dimensional integral  in terms of non-commuting numbers. 
We expect that this direction of research should 
accommodate $q$-calculus on plane and sphere \cite{wess}.

\acknowledgments
We like to thank Professors 
J. Wess, J. Klauder and 
V. Manko  and Dr. K. H. Cho for useful conversations.
This work is supported by the KOSEF
through the CTP at SNU and the project number
(94-1400-04-01-3, 96-0702-04-01-3),
and by the Ministry of Education through the
RIBS (BSRI/96-1419,96-2434).


\begin{references}

\bibitem[*]{poh} E-mail address: ploh@dirac.skku.ac.kr
\bibitem[**]{rim}E-Mail address: rim@phy0.chonbuk.ac.kr
\bibitem{bied} L. C. Biedenharn, J. Phys. A: Math. Gen. {\bf 22}, L873
(1989).
\bibitem{macf} A. J. Macfarlane, J. Phys. A: Math. Gen. {\bf 22}, 4581
(1989).
\bibitem{skly} M. Jimbo, {\it Yang-Baxter equation 
in integrable systems}
 (World Scientific, Singapore, 1989).
\bibitem{Kulda} P. P. Kulish and E. V. Damaskinsky, 
J. Phys. A: Math. Gen. {\bf 23}, L415 (1989).
\bibitem{cho} K. H. Cho, C. Rim, D. S. Soh and S. U. Park,
J. Phys. A: Math. gen. {\bf 27}, 2811 (1994).
\bibitem{hols} T. Holstein and H. Primakoff, Phys. Rev. {\bf 58},
1098 (1940).
\bibitem{dyso} F. J. Dyson, Phys. Rev. {\bf 102}, 1217 (1956).
\bibitem{kitt} C. Kittel, {\it Quantum Theory of Solids},
Second Revised Edition (John Wiley and Sons, New York, 1987).
\bibitem{rowe} A. Klein and E. R. Marshalek, Rev. Mod. Phys. {\bf 63},
375 (1991).
\bibitem{bied1} See L. C. Biedenharn and M. A. Lohe, {\it 
Quantum group symmetry and $q$-tensor algebras}
(World Scientific, Singapore, 1995).
\bibitem{gong} Z. Yu, J. Phys. A: Math. Gen. {\bf 24}, L1321 (1989);
R. Gong, {\it ibid} {\bf 25}, L1145 (1992); D. Ellinas, {\it ibid},
{\bf 26}, L543 (1993).
\bibitem{curt} T. L. Curtright and C. K. Zachos, Phys. Lett. B,
  {\bf 243}, 237.
\bibitem{pere} A. Perelomov, {\it Generalized Coherent States and Their
  Applications} (Springer-Verlag, Berlin, 1986).
\bibitem{klau} J. R. Klauder and B. S. Skagerstam,
{\it Coherent States: Applications in Physics and Mathematical Physics}
 (World Scientific, Singapore, 1985).  
\bibitem{chorim} M. D. Johnson and G. S. Canright, 
Phys. Rev. B {\bf 41}, 6870 (1990);
K. H. Cho and C. Rim, Ann. Phys. (N.Y.), {\bf 213}, 295
(1992).
\bibitem{mac} A. J. Macfarlane, J. Math. Phys. {\bf 35}, 1054 (1994).
\bibitem{gree} H. S. Green,  Phys. Rev. {\bf 90}, 270 (1953);
Y. Ohnuki and S. Kamefuchi, {\it Quantum Field Theory and 
Parastatistics} (Spinger-Verlag, Berlin, 1982).  
\bibitem{oh955} P. Oh, Nucl. Phys. B {\bf 462}, 551 (1996).
\bibitem{rand} S. Randjbar-Daemi, A. Salam and J. Strathdee,
Phys. Rev. B {\bf 48}, 3190 (1993).
\bibitem{sunfu} C. P. Sun and H. C. Fu, J. Phys. A: Math. Gen.
{\bf 22}, L983 (1989); A. Kundu and B. Basu Mallick, Phys. Lett. A
{\bf 156}, 175 91991); J. da Provid\^{e}ncia, J. Phys. A: Math. Gen.
{\bf 26}, 5845 (1993).
\bibitem{baul} L. Baulieu and E. G. Floratos, Phys. Lett. B {\bf 258},
171 (1991); M. Chaichian and A. P. Demichev, Phys. Lett. B {\bf 320},
273 91994).
\bibitem{wess} J. Wess and B. Zumino, Nucl. Phys. Suppl. B {\bf 18},
302 (1990); P. Podl\'{e}s, Lett. Math. Phys. {\bf 18}, 107 (1989).    

\end{references}
\end{document}